\documentclass[a4paper,11pt]{article}
\usepackage{setspace}
\usepackage[english]{babel}
\usepackage{microtype}
\usepackage{amssymb}
\usepackage{mathtools}
\usepackage{times}
\usepackage{amsmath}
\usepackage{bm}
\usepackage{ntheorem}
\usepackage{empheq}
\usepackage{graphicx,psfrag,epsf}
\usepackage{csquotes}
\usepackage{caption}
\usepackage{capt-of}
\usepackage{tabularx}
\usepackage{enumerate}
\usepackage[numbers]{natbib}
\usepackage{url} 
\usepackage[dvipsnames]{xcolor}

\usepackage{lmodern}
\usepackage[utf8]{inputenc}
\usepackage[T1]{fontenc}

\usepackage{graphicx}
\usepackage{xcolor}
\usepackage{tgtermes}
\usepackage{amsfonts}
\usepackage{setspace}
\usepackage{hyperref}
\usepackage[nameinlink]{cleveref}

\usepackage{geometry}
\definecolor{green(html/cssgreen)}{rgb}{0.0, 0.5, 0.0}

\theorembodyfont{\upshape}
\theoremstyle{plain} 
\newtheorem{theorem}{Theorem}[section]

\newtheorem*{Proof}{Proof}
\newtheorem{lemma}[theorem]{Lemma}

\newtheorem{proposition}[theorem]{Proposition}
\newtheorem{corollary}[theorem]{Corollary}

\newtheorem{fact}[theorem]{Fact}

\newtheorem{definition}[theorem]{Definition}
\newtheorem{convendef}[theorem]{Conventional Definition}

\newtheorem{notation}[theorem]{Notation}
\theorembodyfont{\itshape}
\newtheorem{remark}[theorem]{Remark}
\newtheorem{convention}[theorem]{Convention}
\newtheorem{conjecture}[theorem]{Conjecture}

\DeclareMathOperator{\IMM}{Im}
\DeclareMathOperator{\re}{Re}
\DeclareMathOperator{\Arg}{Arg}

\definecolor{blue-violet}{rgb}{0.54, 0.17, 0.89}
\definecolor{ballblue}{rgb}{0.13, 0.67, 0.8}
\definecolor{blue(ryb)}{rgb}{0.01, 0.28, 1.0}
\definecolor{green(html/cssgreen)}{rgb}{0.0, 0.5, 0.0}
\begin{document}
\definecolor{ao(english)}{rgb}{0.0, 0.5, 0.0}
\singlespacing
\title{An exact family of bivariate polynomials and Variants of Chinburg's Conjectures}
\author{Marie-José Bertin \cr Mahya Mehrabdollahei }
\date{} 
\maketitle
\begin{abstract}
This article provides some solutions to 
Chinburg's conjectures by studying a sequence of multivariate polynomials. These conjectures assert that for every odd quadratic Dirichlet Character of conductor $f$, $\chi_{-f}=\left(\frac{-f}{.}\right)$, there exists a bivariate polynomial (or a rational function in the weak version) whose Mahler measure is a rational multiple of $L'(\chi_{-f},-1)$. 
To obtain such solutions for the conjectures we investigate a polynomial family denoted by $P_d(x,y)$, whose Mahler measure has been recently studied.  We demonstrate that the Mahler measure of $P_d$ can be expressed as a linear combination of Dirichlet $L$-functions, which has the potential to generate solutions to Chinburg's conjectures.  Specifically, we prove that this family provides solutions for conductors $f=3,4,8,15,20$, and $24$. 
Notably, $P_d$ polynomials also provide intriguing examples where the Mahler measures are linked to 
$L'(\chi,-1)$ with $\chi$ being an odd non-real primitive Dirichlet character. These examples inspired us to generalize Chinburg's conjectures from real primitive odd Dirichlet characters to all primitive odd characters. For this generalized version of Chinburg's conjecture, $P_d$ polynomials provide solutions for conductors $5,7$, and $9$. 
\end{abstract}
\subsection*{Keywords}
Mahler measure, Polynomial, Chinburg's Conjecture, Dirichlet $L$-function, Dirichlet character, Dilogarithm.
\subsection*{Acknowledgments}
We sincerely thank Antonin Guilloux, Fabrice Rouillier, and François Brunault for their insightful discussions, which enriched the quality of this work. Furthermore, we extend our deepest gratitude to David Hokken, Matilde Lalín, James McKee, Boaz Moerman, Riccardo Pengo, Berend Ringeling, and Wadim Zudilin for their fruitful comments, which enhanced the presentation and clarity of the present paper.
\newpage
\section{Introduction}\label{sec1}
The Mahler measure is a height function that quantifies the complexity of (Laurent) polynomials. 
 \begin{definition}[\normalfont K.~Mahler, \cite{11}]
     Let $P\in \mathbb{C}[z_1,\ldots,z_n]\setminus\{0\}$. The \textit{logarithmic Mahler measure}, $m(P)$, is defined by:
\begin{align}\label{mahler}
    m(P)  \coloneqq  \frac{1}{(2\pi i)^n} \left(\int \cdots \int_{|z_1|=\cdots =|z_n|=1} \log  \bigl |P(z_1, \ldots, z_n) \bigr|  \ \  \frac{dz_1}{z_1}\cdots \frac{dz_n}{z_n} \right),
\end{align}
and the \textit{Mahler measure} of $P$ is defined by $M(P)\coloneqq\exp(m(P))$.
 \end{definition}

Mahler proved that the integral in \Cref{mahler} always exists. For univariate cases, Jensen's equality \cite{Jensen} gives a closed formula for the Mahler measure in terms of the roots of $P$. However, there is no general closed formula for computing the Mahler measure in multivariate cases. Even an error-controlled approximation may be impossible. These
challenging computations are interesting for various reasons such as the links between
Mahler measures and special values of $L$-functions. The second author of the present paper in \cite{mehrabdollahei2021mahler} studied the Mahler measure of a multivariate family of polynomials defined as follows:
\begin{definition}\label{pddef}
    Let $d \in \mathbb{Z}_{\geq 1}$, the polynomial family, denoted by $P_d$, is defined by:
    $$P_d(x,y)\coloneqq\sum_{0 \leq i+j\leq d} x^i y^j.$$
\end{definition}
 The limit of the sequence of the Mahler measure of $P_d$ and its asymptotic behavior is done in \cite{mehrabdollahei2021mahler,BGMP}. The present paper aims to relate individual polynomials of the sequence $(P_d)_{d \in \mathbb{N}}$ to $L$-functions and generate solutions to some important conjectures due to Chinburg namely Conjectures \ref{chinberg}, and \ref{chinberg weak}, as well as a generalization of these conjectures which are respectively \Cref{strong-gen-chin} and \Cref{genchin}. \\

Let $\phi$ be the Euler's totient function and we prove the following theorem, which forms the basis for other results obtained in the present paper.
\begin{theorem}\label{therorm character limit}
Let $d\in \mathbb{Z}_{\geq 1}$, and  for every primitive odd Dirichlet character $\chi$ of conductor $k$, such that $k$ divides $(d+1)(d+2)$, there exists a coefficient $C_{k,\chi}^d\in \mathbb{Q}\left(e^{\frac{2\pi i}{\phi(k)}}\right)\subset \mathbb{Q}\left(e^{\frac{2\pi i}{\phi((d+1)(d+2))}}\right)$ such that:
$$m(P_d)=\sum_{k|(d+1)(d+2)}\ \sum_{\chi \text{ primitive odd mod } k} C_{k,\chi}^d L'(\chi,-1).$$
\end{theorem}
\Cref{tablemp2} illustrates the above theorem for $m(P_d)$ with $d\leq 16$. For the constants $C_{k,n}^d$, associated with the coefficient of $L'(\chi_{k}(n,.),-1)$ in $m(P_d)$, appearing in the following table, we refer to \Cref{table-of-constant} in \Cref{appendix}. The characters that appeared in \Cref{table-of-constant} are expressed using the Conrey representation in the LMFDB, \cite[Section 2]{Conrey} ,which is also recalled in \Cref{appendix}. 
\begin{center}

\fbox{%
    \parbox{\textwidth}{%
\small{\begin{align*}
&m(P_1)=L'(\chi_{-3},-1)\\
&m(P_2)=L'(\chi_{-4},-1)-\frac{L'(\chi_{-3},-1)}{2}\\
&m(P_3)=\re\left(\frac{9-3i}{10}L'(\chi_5(2,.),-1)\right)-\frac{3}{5}L'(\chi_{-4},-1)\\
&m(P_4)=\re\left(\frac{-3+i}{5}L'(\chi_5(2,.),-1)\right)+\frac{16}{5}L'(\chi_{-3},-1)\\
&m(P_5)=\frac{1}{3}L'(\chi_{-7},-1) +\re\left( \frac{8-4\sqrt{3}i}{21}L'(\chi_{7}(3,.),-1)\right) -\frac{16}{7}L'(\chi_{-3},-1)\\
&m(P_6)=\frac{3}{7}L'(\chi_{-8},-1)-\frac{1}{4}L'(\chi_{-7},-1)+\re\left(\frac{-2+\sqrt3i  }{7}L'(\chi_{7}(3,.),-1)\right)+\frac{15}{14}L'(\chi_{-4},-1) \\
&m(P_7)=\re\left(\frac{3-\sqrt{3}i}{6}L'(\chi_9(2,.),-1)\right)-\frac{1}{3}L'(\chi_{-8},-1)-\frac{5}{6}L'(\chi_{-4},-1)+\frac{5}{6}L'(\chi_{-3},-1)\\
&m(P_8)=\re\left(\frac{-6+2\sqrt{3}i}{15}L'(\chi_9(2,.),-1)\right)+\re\left(\frac{28-16i}{15}L'(\chi_5(2,.),-1)\right)-\frac{2}{3}L'(\chi_{-3},-1)\\
    &m(P_9)=\frac{1}{5}\re\left(C_{11,2}^9L'(\chi_{11}(2,.),-1)+C_{11,7}^9L'(\chi_{11}(7,.),-1)-\frac{84-48i}{11}L'({\chi_5(2,.)},-1)\right)+\frac{3}{25}L'(\chi_{-11},-1)\\
&m(P_{10})=-\frac{1}{6}\re\left(C_{11,2}^{10}L'(\chi_{11}(2,.),-1)+C_{11,7}^{10}L'(\chi_{11}(7,.),-1)\right)\\
    & \qquad\quad\quad-\frac{1}{10}L'(\chi_{-11},-1)+\frac{1}{11}\left(21L'(\chi_{-4},-1)+38L'(\chi_{-3},-1)\right)\\
 &m(P_{11})=\frac{1}{6}\re\left(C_{13,5}^{11}L'({\chi_{13}(5,.)},-1)+C_{13,2}^{11}L'(\chi_{13}(2,.),-1)+C_{13,6}^{11}L'(\chi_{13}(6,.),-1)\right)\\
    & \qquad\quad\quad-\frac{1}{13}\left(21L'(\chi_{-4},-1)+38L'(\chi_{-3},-1)\right)\\
    &m(P_{12})=-\frac{1}{7}\re\left(C_{13,5}^{12}L'({\chi_{13}(5,.)},-1)+C_{13,2}^{12}L'(\chi_{13}(2,.),-1)+C_{13,6}^{12}L'(\chi_{13}(6,.),-1)\right)
\\& \quad\quad\quad\quad +\re\left(\frac{92-60 \sqrt{3}i}{91} L'({\chi_{7}(3,.)},-1)\right)+ \frac{4}{13}L'(\chi_{-7},-1)\\
&m(P_{13})=\frac{3}{14}L'({\chi_{-15}},-1)+\re\left(\frac{-92+60 \sqrt{3}i}{105} L'({\chi_{7}(3,.)},-1)\right)
\\& \quad\quad\quad\quad 
    -\frac{4}{15}L'({\chi_{-7}},-1)+\re\left(\frac{48-6i}{35}L'({\chi_5(2,.)},-1)\right)+\frac{8}{7}L'({\chi_{-3}},-1)\\
&m(P_{14})=\re\left(\frac{1+i}{5}L'(\chi_{16}(3,.),-1)\right)-\frac{3}{16}L'(\chi_{-15},-1)+\frac{1}{2}L'(\chi_{-8},-1)\\
        &\qquad\qquad  -\re\left(\frac{24-3i}{20}L'({\chi_5(2,.)},-1)\right)+\frac{21}{20}L'(\chi_{-4},-1)-L'(\chi_{-3},-1).\\
&m(P_{15})=\frac{1}{8}\re\left(C_{17,3}^{15}L'({\chi_{17}(3,.)},-1)+C_{17,5}^{15}L'({\chi_{17}(5,.)},-1)+C_{17,10}^{15}L'({\chi_{17}(10,.)},-1)+C_{17,12}^{15}L'({\chi_{17}(12,.)},-1)\right)\\
        &\qquad\qquad -\re\left(\frac{3+3i}{17}L'(\chi_{16}(3,.),-1)\right)-\frac{15}{34}L'(\chi_{-8},-1)-\frac{63}{68}L'(\chi_{-4},-1).\\
        &m(P_{16})=
        -\frac{1}{9}\re\left(C_{17,3}^{16}L'({\chi_{17}(3,.)},-1)+C_{17,5}^{16}L'({\chi_{17}(5,.)},-1)+C_{17,10}^{16}L'({\chi_{17}(10,.)},-1)+C_{17,12}^{16}L'({\chi_{17}(12,.)},-1)\right)
\\&\qquad\qquad 
        +\re\left(\frac{36-20\sqrt{3}i}{51}L'(\chi_9(2,.),-1)\right)+\frac{160}{51}L'(\chi_{-3},-1)
            \end{align*}}}}

\captionsetup{type=figure}\captionof{table}{Table of the representation of $m(P_d)$ in terms of $L$-functions, for $1\leq d \leq 16$.}
\label{tablemp2}
\end{center}
\newpage
The present paper explores Chinburg's conjectures on the links between Mahler measures and special values of Dirichlet $L$-functions. We note that all versions of Chinburg's conjecture should be considered alongside \Cref{negative odd integeres}.

\begin{conjecture}[\normalfont Strong Chinburg's Conjecture \cite{chinburge}]\label{chinberg}
For every odd quadratic character $\chi_{-f}\coloneqq\left(\frac{-f}{.}\right)$, and $n \in \mathbb{N}$ there exists a non-zero polynomial $Q_f\in \mathbb{Z}[z_1,\ldots, z_{n+1}]$ and a rational number $r_f$ for which $$r_fm(Q_f)=L'(\chi_{-f},-n). $$ 
\end{conjecture}
We recall that $\left(\frac{-f}{.}\right)$ is the Kronecker symbol. Chinburg \cite{chinburge} mentioned a weak version of the above conjecture as well.
\begin{conjecture}[\normalfont Weak Chinburg's Conjecture \cite{chinburge}]\label{chinberg weak}
For every odd quadratic character $\chi_{-f}$, there exists a non-zero rational function $R_f\in \mathbb{Q}(z_1,\ldots, z_{n+1})$ and a rational number $r_f$ for which $$r_fm(R_f)=L'(\chi_{-f},-n). $$
\end{conjecture}
In \cite{chinburge} there is proof for the weak version, for the case $n=1$, but the proof appears to be incorrect, therefore both weak and strong versions of Chinburg's problem remain open.
\Cref{therorm character limit} expresses the Mahler measures
$m(P_d)$ in terms of special values of $L$-functions. Conjectures \ref{chinberg}, and \ref{chinberg weak}) aim to express special values of $L$-functions as linear combinations of Mahler measures. In the present paper, using the expression of $m(P_d)$ in terms of $L$-functions, we construct
sequences $a = (a_d)_d \in \mathbb{Z}^n$ with finite support such that the Mahler measure of the rational 
function $\Pi_{d=0}^{\infty}P_d^{a_d}$ is a rational multiple of $L'(\chi_{-f},-1)$, where $\chi_{-f}$ is the odd quadratic Dirichlet character of conductor $f$. More precisely, we construct the following partial solutions to Chinburg's conjectures;

\renewcommand{\arraystretch}{3}

\begin{center}
\fbox{%
    \parbox{\textwidth}{%
\small{
\begin{align}
&m(P_1)=L'(\chi_{-3},-1)\label{d3}\\
&\nonumber\\
&m(P_1P_2^2)=2L'(\chi_{-4},-1)\label{d4} \\
&\nonumber\\
&m(\dfrac{P_1^{33}\,P_5^{21}\,P_6^{28}}{P_2^{30}})=12 L'(\chi_{-8},-1)\label{d8}\\
&\nonumber\\
&m(\dfrac{P_{13}^{210}\,P_{12}^{182}\,P_{11}^{156}\,P_{8}^{135}\,P_7^{108}\,P_6^{84}\,P_5^{63}\,P_4^{900}\,P_2^{252}}{P_{1}^{2394}})=45 L'(\chi_{-15},-1)\label{d15}\\
&\nonumber\\
&m\left(\dfrac{P_{18}^{190}\,P_{17}^{171}\,P_{16}^{153}\,P_{15}^{136}\,P_{14}^{120}\,P_{13}^{105}\,P_{12}^{91}\,P_{11}^{78}\,P_{1}^{168}}{P_{8}^{225}\,P_{7}^{180}\,P_{6}^{140}\,P_{5}^{105}\,P_{2}^{24}\,P_{4}^{60}}\right)=30L'(\chi_{-20},-1)\label{d20}\\
&\nonumber\\
&m\left(\dfrac{P_{22}^{276}\,P_{21}^{253}\,P_{20}^{231}\,P_{19}^{210}\,P_{18}^{190}\,P_{17}^{171}\,P_{16}^{153}\,P_{15}^{136}\,P_{14}^{120}\,P_{13}^{105}\,P_{12}^{91}\,P_{11}^{78}}{P_{6}^{252}\,P_{5}^{189}\,P_{2}^{234}\,P_{1}^{1269}}\right)=36L'(\chi_{-24},-1)\label{d24}
\end{align}
}}}
\captionsetup{type=figure}\captionof{table}{Table of some solutions to Chinburg's conjectures using the $P_d$\, polynomials.}\label{6.3.6}
\end{center}
The case $f=3$ in \Cref{6.3.6} was first done by Smyth \cite{smyth_1981}, but we recover this result using $P_d$ family. Furthermore, $P_d$ family provides interesting examples involving the $L$-functions of complex primitive characters that suggest generalizing Chinburg's conjecture to all primitive odd characters as follows;
 \begin{conjecture}[\normalfont Generalization of the weak Chinburg's Conjecture]\label{genchin}
    For every primitive odd Dirichlet character $\chi$ of conductor $f$, there exists a non-zero rational function $R_\chi \in \mathbb{Q}(z_1,\ldots, z_{n+1})$ and a rational number $r_\chi\ $ for which
    \begin{align}\label{chin-gen-eq}
        r_\chi m(R_{\chi})=2 \re(L'(\chi,-n)).
    \end{align}
    
\end{conjecture}
\begin{conjecture}[\normalfont Generalization of the strong Chinburg's Conjecture]\label{strong-gen-chin}
   For every primitive odd Dirichlet character $\chi$ of conductor $f$, there exists a non-zero polynomial $Q_\chi \in \mathbb{Z}[z_1,\ldots, z_{n+1}]$ and a rational number $r_\chi\ $ for which
    \begin{align}\label{chin-gen-eq-strong}
        r_\chi m(Q_\chi)=2 \re(L'(\chi,-n)).
    \end{align}
     
\end{conjecture}
We note that the Mahler measure is a real number, so generalizing Chinburg's conjectures to complex characters obliged us to involve the $L$-functions associated with $\chi$ and its conjugate $\bar{\chi}$. Moreover, if one can construct $R_\chi$ for a given Dirichlet character that satisfies \Cref{chin-gen-eq}, then one can take $R_\chi=R_{\bar{\chi}}$.  If in the above conjecture, we let $\chi$ be a real primitive character (i.e. quadratic), then $\bar{\chi}=\chi$ and we recover the actual version of Chinburg's conjectures.
Using the family $\{P_d\}_d$ we could only provide solutions to the weak version of the generalization of Chinburg's conjecture.
\begin{proposition}
     Let $Q$ be the polynomial  introduced by Ray \cite{Ray1987RelationsBM} 
as a solution to the strong Chinburg's conjecture for conductor $7$, satisfying $\frac{8}{7}L'(\chi_{-7},-1)=m(Q)$,  as follows:\begin{align}\label{Ray-7-poly}
  Q(x,y)=&\frac{(x^7-1)}{x-1}(y-1)^2+7x^2(x+1)^2y  \end{align}
\Cref{genchinsolution} lists the solutions to \Cref{genchin}, generated mainly by the $P_d$ family.

\renewcommand{\arraystretch}{2.2}
\begin{center}
\fbox{%
    \parbox{\textwidth}{%
\small{
\begin{align}
&m\left(\dfrac{P_1^{720}}{P_8^{45}\,P_7^{36}\,P_6^{28}\,P_5^{21}\,P_4^{240}}\right)=60\re\left(L'(\chi_5(2,.),-1)\right)\label{dchi5}\\
    &\nonumber\\
&m\left(\dfrac{P_1^{3432}\,P_5^{2520}}{P_{12}^{728}\,P_{11}^{624}\,P_2^{1008}\,Q^{539}}\right)=224\re\left(L'(\chi_7(3,.),-1)\right)\label{dchi7}\\
&\nonumber\\
&m\left(\dfrac{P_{7}^{360}\,P_{6}^{280}\,P_{5}^{210}\,P_{1}^{369}}{P_{16}^{153}\,P_{15}^{136}\,P_{14}^{120}\,P_{13}^{105}\,P_{12}^{91}\,P_{11}^{78}\,P_{2}^{126}}\right)=72\re\left(L'(\chi_9(2,.),-1)\right)\label{dchi9}
\end{align}}}}
\captionsetup{type=figure}\captionof{table}{Table of solutions to \Cref{genchin} using the $P_d$\, family.}
\label{genchinsolution}
\end{center}
\end{proposition}

\begin{remark}
    All the identities in \Cref{tablemp2}, \Cref{6.3.6} and \Cref{genchinsolution} have been numerically verified. The SageMath code we used to verify the identities is available online
    \cite{code}.
\end{remark}

\begin{remark}\label{negative odd integeres}
The present paper restricts its scope to the case of two-variable polynomials and, in the context of Chinburg's conjectures, focuses on evaluating the derivatives of Dirichlet $L$-functions at $-1$. In the general case, it is worth noting that the oddness of the character $\chi$ implies that $L(\chi, s)$ has trivial zeros at negative odd integers. These zeros are
all simple, which implies that for every $n \geq 1$ the derivative $L'(\chi,1-2n)$ coincides with the
special value $L^*(\chi,1-2n)$, where $L^*(\chi,1-2n)$ is the completed Dirichlet $L$-function (see for instance \cite{completed-L}). On the other hand, if $\chi$ is a primitive odd Dirichlet character we
have that $L^*(\chi,1-2n)=L(\chi,-2n)=\frac{-\beta_{k,\chi}}{k} \in \mathbb{Q}(\chi)^\times$,
where $\{\beta_{k,\chi}\}_k$ denotes the sequence of
generalized Bernoulli numbers associated with $\chi$.
If $\chi$ is a
primitive odd Dirichlet character the number $L'(\chi,-2n)$ could even be expected not to be a
period, in the sense of Kontsevich-Zagier. As a consequence,  we restrict all versions of Chinburg's conjectures to negative odd integers. 
\end{remark}
\begin{remark}\label{mostgeneralchin}
    There is a possibility of further generalizing Chinburg's conjecture to include all Dirichlet characters. However, this broader generalization is likely to involve algebraic coefficients rather than rational ones. We believe that constructing solutions to the most general version of Chinburg's conjecture using the 
$P_d$ family and polynomials suggested by Ray in \cite[Proposition 22]{Ray1987RelationsBM}, could provide valuable insights. However, this lies beyond the scope of the present paper.

\end{remark}
 \subsection{Historical remarks}

 We devote this subsection to a short historical overview on  the existing
results regarding Chinburg's conjectures.\\

 Historically, the first  explicit
formulas linking Mahler measures with special values were obtained by C.J.~Smyth \cite{smyth_1981} \begin{equation}\label{firstsmyth}
  m(1+x+y)=\frac{3\sqrt{3}}{4\pi}L(\chi_{-3},2)=L'(\chi_{-3},-1),
    \end{equation}
 and  
 \begin{equation}\label{secondsmyth}
 m(1+x+y+z)=\frac{7}{2\pi^2}\zeta(3)=-14\zeta'(-2),
  \end{equation}
 where $\chi_{-3}(n)=\left(\frac{-3}{n}\right)$ is the odd quadratic Dirichlet character of conductor $3$ and $\zeta$ is the Riemann Zeta function. Inspired by \Cref{firstsmyth}, Chinburg introduced the conjecture formalized in \Cref{chinberg}. 
For further fascinating insights into multivariate Mahler measure computations and connections to special values we direct the reader to the works of Bertin \cite{5,unemesure,mahlerbertin,MR2076562}, Boyd \cite{7,9,Boydinvariant}, Boyd and Rodriguez-Villegas \cite{Mahler2,boyd_rodriguez-villegas_2002}, Rodriguez-Villegas \cite{Villegas1999}, Smyth \cite{smyth_1981b,Anexplicitformula}, Bertin and Zudilin \cite{4,3}, and Lalín \cite{lalinvariable,LalinSomeexamples,Lalincombination,8,6}.\\
A fundamental connection between Mahler measures and special values of $L$-functions arises through the Bloch-Wigner dilogarithm (see \Cref{def 2.3}), as the Mahler measure of certain polynomials often involves sums of the dilogarithm at roots of unity, which are directly connected to special values of Dirichlet $L$-functions.

\begin{proposition}\label{df}{\cite{bloch2000higher}}
Let $-f$ be a fundamental discriminant and $\chi_{-f}\coloneqq\left( \frac{-f}{.}\right)$ be the odd quadratic Dirichlet character of conductor $f$. Boyd, \cite{7}, introduced the notation $d_f$ for which the following equalities are valid:
\begin{align}\label{notation df}
d_{f}\coloneqq\frac{f}{4\pi}\sum_{m=1}^{f-1}\chi_{-f} (m)D(e^{\frac{2\pi i m}{f}})=\frac{f^{\frac{3}{2}}}{4\pi}L(\chi_{-f},2)= L'(\chi_{-f},-1), 
\end{align}
where $D(z)$ is the Bloch-Wigner dilogarithm function  which is introduced in \Cref{def 2.3}.
\end{proposition}

Chinburg's conjectures have been explored exclusively for the case $n=1$ and resolved for just 18 specific conductor values. In the present paper, we also focus solely on the case $n=1$.\\
Ray \cite{Ray1987RelationsBM} constructed polynomials answering to \Cref{chinberg}, for \( f=3, 4, 7, 8, 20, \) and \( 24 \). \Cref{Table of Rayexamples} in \Cref{appendix} shows the polynomials constructed by Ray and the explicit value of the constants $r_f$. His proof for \( f=7 \)  is more intricate, relying on linear relations among certain 
 kinds of twisted 
$L$-series. For the precise definition of these twisted 
$L$-series, see \cite[Page 698]{Ray1987RelationsBM}. The method he applied doesn't seem to work for other characters. In  \cite[Proposition 22]{Ray1987RelationsBM} Ray establishes connections between 
$L$-functions of certain odd Dirichlet characters with even order and Mahler measures. Notably, his results could apply to primitive Dirichlet characters 
 presented in \Cref{genchinsolution}. However, by applying Ray's method we obtain polynomials whose Mahler measures are of the form $\re(C_\chi L'(\chi,-1))$, where $C_\chi$ is an algebraic number and not necessarily rational. In particular for $\chi=\chi_5(2,.),\ 
\chi_7(3,.)$ and $\chi_9(2,.)$ we have $C_\chi=\frac{8+2i}{10},\frac{8-2e^\frac{4 \pi i}{3}}{21}$ and $\frac{16-4e^\frac{5 \pi i}{3}}{54}$, respectively. In contrast to Ray's work, all the polynomials that we constructed in \Cref{genchinsolution} have $C_\chi\in \mathbb{Q}$. \\

Furthermore, another solution to the strong Chinburg's conjecture for conductor $7$ with an explicit constant is provided by the first author of the present paper in \cite{Bertinlink1991} as follows:
$$d_7=3m((x+1)^2 y+x^2+x+1).$$
The proof of the above equality is an application of Vandervelde's general formula \cite{VANDERVELDE20082231} for the Mahler measure of a polynomial $P(x,y)$ defining a curve of genus zero, combined with the dilogarithmic identity $D(a^3) = 3D(a)-D(a/2)$ involving the quadratic irrational $a$ with minimal polynomial $x^2 + \frac{3}{2} x + 1$.
There is also some of her unpublished
work regarding conductor $15$ that is announced in conferences in \href{https://webusers.imj-prg.fr/%7Emarie-jose.bertin/edmonton08.pdf}{2008} and \href{https://webusers.imj-prg.fr/%7Emarie-jose.bertin/montreal15.pdf}{2015}, as follows:
\begin{align*}
&12m\left((x^2+x+1)y+x^2+1\right)=d_{15},\\
    & m\left((x^2+x+1)(y^2+x^2)-(x^4-x^3-6 x^2-x+1)y\right)=\frac{1}{3}d_7+\frac{1}{6}d_{15}.
\end{align*}
In the proof of the above equalities, she mainly uses Vandervelde's formula \cite{VANDERVELDE20082231}, the Bloch group of the number field of toric points, and the expression of the zeta function of a number field in terms of the Bloch-Wigner dilogarithms. In particular, for the first equation, she applied the dilogarithmic identity $D(b^3)-3D(b)-2D(-b)+2D(b/2)=0$, where $b$ has the minimal
polynomial $2x^2+x+2$. \\

Boyd and Rodriguez-Villegas, in \cite{boyd_rodriguez-villegas_2002} and \cite{Mahler2} provided solutions to the strong Chinburg's conjecture, for $f= 3, 4, 7, 8, 11, 15, 19, 20, 24,$ $ 35, 39, 40, 55, 84,$ and $ 120$. Their method involves proving that the Mahler measure of a specific family of polynomials of the form $p(x)y-q(x)$, where $p(x)$ and $q(x)$ are cyclotomic polynomials, is the Bloch-Wigner dilogarithm of an element of the Bloch group of a certain quadratic imaginary field. Then, due to a rank restriction, the Mahler measures in question must be a multiple of certain $d_f$, where $f$ is related to the polynomial's discriminant. The rational constant $r_f$ in their method is approximated using computer programs. For this reason, we use the notation \( d_f \sim r_f m(P_f) \) in \Cref{Numericalyverifiedtable}. One important polynomial family they fully investigate, which is obtained by setting $p(x)=1$ in $p(x)y-q(x)$, is the genus-zero family of the form 
 $y+1+x+\ldots+x^d$, for $d\in \mathbb{N}$, where the degree is unbounded. This family is exact (see \Cref{Exactdef}) and in \cite[Proposition 1]{boyd_rodriguez-villegas_2002} they compute a closed formula for the Mahler measure of this family. They indeed show certain relations between  $d_f$'s and the Mahler measure of member of this polynomial family. Liu and Qin in \cite{Luo_2019} studied more than 500 families of reciprocal polynomials defining genus 2 and genus 3 curves. They found numerical relations between the Mahler measures of these polynomials and special values of $L$-functions. They
gave solutions to the strong Chinburg's conjecture for $f= 23, 303,$ and $755$. 
\Cref{Numericalyverifiedtable} in \Cref{appendix} summarizes all the examples with numerically approximated constants that we have thanks to \cite{boyd_rodriguez-villegas_2002,Mahler2,Luo_2019}.


\section{The Mahler measure of \texorpdfstring{$P_d$}{} family and links with \texorpdfstring{$L$}{}-functions}\label{se-2}
As we explained in the introduction, we aim to provide solutions to Chinburg's conjectures using the polynomial family $P_d$. Part of the work of \Cref{se-2} and \Cref{Solutions to Chinburg's conjecture-sec} is done in the doctoral thesis of the second author, \cite{phdthesis}. We first demonstrate that $m(P_d)$ can be expressed as linear combinations of $L$-functions. This is the main result of this chapter, stated as \Cref{therorm character limit} in the introduction.
Before delving into the proof of the theorem, we provide a concise introduction to the $P_d$ family, emphasizing its most significant properties and distinctive features that already exist in the literature.
\subsection{A brief introduction to \texorpdfstring{$P_d$}{} 
 family}
 
The polynomial family $P_d$, introduced in \Cref{pddef}, can be expressed in the following summation-free form:
\begin{align*}
P_d(x,y)=\frac{(y-1)x^{d+2}-(x-1)y^{d+2}+(x-y)}{(1-x)(1-y)(x-y)}.
\end{align*}
 The family $(P_d)_{d \in \mathbb{N}}$ is projectively smooth and unbounded in the total degree (both variables $x$ and $y$). The genus-degree formula implies that the genus of the algebraic curve $P_d$ is $\frac{(d-1)(d-2)}{2}$. An alternative method to compute the genus of $P_d$ is explained in  \cite[Proposition VII.5.11]{phdthesis}. Therefore the family $(P_d)_{d \in \mathbb{N}}$ is indeed unbounded in genus.
The Newton polytope of $P_d$ is the triangle with vertices $\{(0,0),(0,d),(d,0)\}$ (See \cite[Example VII.5.6]{phdthesis}).
\begin{definition}[\normalfont \cite{Villegas1999}] Let $P(x,y)= \sum_{(m,n) \in \mathbb{Z}^2}c_{(m,n)}x^my^n \in \mathbb{C}[x^{\pm 1},y^{\pm 1}]$.
A side of $N_P$, Newton polytope of $P$,  is denoted by $\tau$. We parameterize a side of 
$N_P$ clockwise, labeling the consecutive lattice points on $\tau$
as $\tau_{(0)}, \tau_{(1)}, \ldots$.
To each side, we associate a one-variable polynomial, called the \textit{side polynomial}, denoted by $P_{\tau}(t)$, defined as
\begin{align*}
    P_{\tau}(t):=\sum_{k=0}^{\infty} c_{\tau(k)}t^k \in \mathbb{C}[t] 
\end{align*}
This sum is naturally finite. A polynomial whose Newton polytope contains side polynomials with only roots of unity is called \textit{tempered}.
\end{definition}

The polynomial associated with each side of $N_{P_d}$ is $1+t+t^2+ \cdots +t^d$. Therefore, $P_d$, for each $d\in \mathbb{N}$ is a tempered polynomial. 
\begin{convention}
In this article, $\log z$ denotes the natural logarithm, $\IMM(z)$, and $\re(z)$ are respectively the real part and imaginary part of the complex number $z$. Furthermore, the notation $\Arg$ is used for the principal value of the argument that lies within the interval $(-\pi,\pi]$.
 \end{convention}
In the following, we briefly recall the notion of bivariate exactness.   
\begin{definition}[\normalfont \cite{1}]\label{Exactdef}
    A polynomial $P\in \mathbb{C} [x,y]$ is called \textit{exact} if the form $$\eta\coloneqq \log|y|\ d\arg(x)-\log |x|\ d\arg(y)$$ restricted to the algebraic curve $$C_P= \{(x,y)\in {{\mathbb{C}}^*}^2| P(x,y)=0 , dP(x,y) \neq 0 \}$$  is exact. 
\end{definition}

In \cite[Theorem 3.9]{mehrabdollahei2021mahler} is proved that, for every $d \in \mathbb{N}$, the polynomial $P_d$ is exact. 
To present the explicit formula for $m(P_d)$, we need to recall the definition of the standard dilogarithm function as well as the Bloch-Wigner dilogarithm function.
\begin{definition}[\normalfont \cite{2}]\label{stdilog}
    The dilogarithm function is the function defined by the power series
    \begin{equation*}
    \operatorname{Li_2(z)} = \sum_{n=1}^{\infty}\frac{z^n}{n^2}\ , \ \text{for} |z|<1.
\end{equation*}
\end{definition}
The definition and name originate from the analogy with the Taylor series expansion of the ordinary logarithm around $1$.  The analytic continuation of
the dilogarithm is given by the following integral:
\begin{align}\label{eqLI2}
    \operatorname{Li_2(z)} = - \int_{0}^{z} \log(1-u)\frac{du}{u} \ \ \  \ \text{for} \ z \in \mathbb{C} \setminus [1,\infty).
\end{align} 

\begin{definition}\label{def 2.3}
 The \textbf{Bloch-Wigner} dilogarithm $D(z)$ is defined by: 
 
\begin{equation*}
    \operatorname{D(z)} = \IMM(Li_2(z)) + \Arg(1-z)\log |z|.
\end{equation*}
\end{definition}
In the following facts, we summarize the important properties of the Bloch-Wigner dilogarithm. For more details, we refer to \cite{2}.
\begin{fact}\label{dilogarithm}
 The function $D(z)$ is real analytic on $\mathbb{C}$ except at the two points $0$ and $1$, where it is continuous but not differentiable. Moreover, we have:
\begin{enumerate}
    \item $D(\bar{z})=-D(z)$.
    \item If $\theta \in \mathbb{R}$ and $z=e^{i\theta}$ then $D(z)=D(e^{i\theta})=Cl_2(\theta)$, where $$Cl_2(\theta):=-\int_{0}^{\theta}\log| 2\sin(\theta/2)|$$
denotes the Clausen function of order two, which admits the Fourier series representation $$Cl_2(\theta)=\sum_{n=1}^{\infty}\frac{\sin (n\theta)}{n^2}.$$
In particular, for every $k \in \mathbb{Z}$, we have $D\left(e^{k\pi i}\right)=0$.
\item Let $n\in \mathbb{Z}_{\geq 1}$, the following relation holds: $$D\left(z^n\right)=n \sum_{k=0}^{n-1} D\left(e^\frac{2k\pi i}{n} z\right).$$ 
\end{enumerate}
The last property above is called the distribution relation for dilogarithm and is mentioned in \cite[Page 9]{2}.
\end{fact}

\begin{convention}
 In the present paper, we exclusively address the Bloch-Wigner dilogarithm. Therefore, any reference to the term "dilogarithm" should be understood as pertaining specifically to the Bloch-Wigner dilogarithm.  
\end{convention}
In \cite{1}, the authors derived a closed formula for computing the Mahler measure of a class of bivariate exact polynomials. In \cite[Proposition 4.9]{mehrabdollahei2021mahler},  the second-named author obtained a closed formula for the Mahler measure of $P_d$ by utilizing the Mahler measure formula for exact polynomials, yielding the following result:
\begin{align}\label{eq 4.1.11}
      2\pi m(P_d)=\frac{1}{d+1}S_{d+2}-\frac{1}{d+2}S_{d+1}, \quad \text{with}\quad S_d\coloneqq3\sum_{1 \leq k \leq d-1}(d-2k)D\left(e^{\frac{2\pi i k}{d}}\right). 
\end{align}
The $P_d$ family is the first family with unbounded degree and genus for which a closed formula for Mahler measures exists.
By setting $d=1$ we recover the polynomial $P_1=x+y+1$, the well-known polynomial studied by Smyth and introduced in \Cref{firstsmyth}, whose Mahler measure is a special value of a Dirichlet $L$-function. Thus, one can consider $P_d$ as a generalization of $P_1$. \\

This family is indeed fruitfully regarding one of the important conjectures due to Boyd;
\begin{conjecture}[\normalfont Boyd's Conjecture \cite{9}]\label{boydcon} $L^\sharp\coloneqq\{M(P)| P\in \mathbb{Z}[z_1,\ldots, z_n]\setminus\{0\}, \ n\geq 1 \}$ is closed with respect to the Euclidean topology.
\end{conjecture}
Boyd \cite{9} investigated sequences of multivariate Mahler measures and their limits to address Lehmer's famous conjecture, \cite{Lehmer}, regarding a universal (positive) lower bound for the logarithmic Mahler measures of integral polynomials with non-vanishing Mahler measures. Although this approach did not lead to a proof of Lehmer's conjecture, Boyd proposed the conjecture above, which, if proven true, would imply Lehmer's conjecture. However, the computation of the Mahler measure of multivariate polynomials is challenging, and as a result, there is limited evidence for Boyd's conjecture. One of the instances where Boyd's conjecture is seen to hold is for the sequence of  Mahler measures of $P_d$ polynomials. From \cite{mehrabdollahei2021mahler} and the computation of the Mahler measure of the polynomial \begin{align}\label{Pinfty}
P_\infty \coloneqq (1-x)(1-y)-(1-z)(1-w),    
\end{align}

done by D'Andrea
and Lalín in \cite{6},  we have 
\begin{align}\label{limitpd}
    \lim_{d \rightarrow \infty}m(P_d)= m(P_\infty)= \frac{9}{2 \pi^2}\zeta(3).
\end{align} 
This provides an example of a Cauchy sequence of Mahler measures of integral polynomials that converges to a Mahler measure of a four-variable integral polynomial, supporting Boyd's conjecture.\\  

 A straightforward proof of \Cref{limitpd} can be obtained by applying the generalized Boyd-Lawton theorem, as fully explained in \cite{BGMP}. 
The convergence of the sequence of $m(P_d)$ to $m(P_\infty)$ is of the order $O\left(\frac{\log d}{d^2}\right)$, but for the full asymptotic expansion of $|m(P_d) - m(P_\infty)|$ we refer to \cite[Theorem 5.1]{BGMP} or \cite[Theorem V.4.1]{phdthesis}.
\\

   Lastly, we would like to make the following remark regarding $P_d$ family.
\begin{remark}
 The set  $\{m(P_d)|d \in \mathbb{N}\}$ contains $\mathbb{Q}$-linearly dependent elements. For instance, an easy computation reveals the following identities:
 \begin{align*}
  &m(P_4) = \frac{-2}{3}m(P_3) - \frac{2}{5}m(P_2) + 3m(P_1),\\
  &66m(P_{10})+55m(P_{9})+45m(P_{8})+36m(P_{7})+28m(P_{6})+21m(P_{5})=243m(P_{1})+126m(P_{2})
 \end{align*} 
 \end{remark}
\subsection{Representation of \texorpdfstring{$m(P_d)$}{} in terms of \texorpdfstring{$L$}{}-functions}\label{representationpd-sec}
The closed formula for the Mahler measure of the $P_d$, \Cref{eq 4.1.11}, is a linear combination of the Bloch-Wigner dilogarithm function evaluated at roots of unity. 
Furthermore, there is a link between $L$-functions of odd quadratic characters and Bloch-Wigner dilogarithms evaluated at roots of unity as mentioned in \Cref{df}. Based on Smyth's computation, \Cref{firstsmyth}, and \Cref{df}, one may ask about the possibility of writing $m(P_d)$, for $d\geq 2$ as a  a linear combination of  $L'(-1,\chi_{-f})$ for some values of $f$, where $\chi_{-f}$ is the odd quadratic character of conductor $f$. In this section, we prove that $m(P_d)$ can be written as a linear combination of $L$-functions of primitive odd characters (not necessarily quadratic). This result was stated as \Cref{therorm character limit} in the introduction.
In the rest of this section, we provide the proof of the theorem.  We begin by recalling the following definition:
\begin{definition}
Let $\chi$ be a primitive character modulo $k$. The Gauss sum associated to $\chi$ is defined as follows:
 \begin{align*}
     \tau(\chi)=\sum_{1\leq a \leq k} \chi(a) e^{\frac{2\pi i a}{k}}.
 \end{align*}
 \end{definition}
 It is a classical fact  \cite[Page 84]{lang2013algebraic} that for the primitive Dirichlet character $\chi$ of conductor $k$ we have 
$$|\tau(\chi)|^2=k.$$  
The first step in finding a link between $m(P_d)$ and $L$-functions is the following proposition.
 \begin{proposition}\label{pro 5.1}
\cite[Section 9.B]{Zagier1991} Let $\chi$ be a complex primitive odd Dirichlet character of conductor $k$. We have:
\begin{align}\label{L(2)}
    L(\chi,2)=i \tau(\bar{\chi})^{-1}\sum_{m=1}^{ k-1} \overline{\chi(m)}D\left(e^{\frac{2\pi i m}{k}}\right).
\end{align}
\end{proposition}
We recall that $L$ -functions satisfy the following functional equation for an odd primitive Dirichlet character $\chi$;

\begin{align}\label{feL}
  L'(\chi,-1)=  \frac{ -ik\tau(\chi)}{4\pi} L(\bar{\chi},2)
\end{align}
In the following corollary, using \Cref{feL} and \Cref{pro 5.1}, one can generalize the \Cref{df} to all primitive Dirichlet characters.
\begin{corollary}[\normalfont \cite{Ray1987RelationsBM}, Page 697]
\label{coclusion1.2}
For a primitive odd (non-principal) Dirichlet character of conductor $k$, we have:
$$L'(\chi,-1)=\frac{ -ik\tau(\chi)}{4\pi} L(\bar{\chi},2)=\frac{k}{4\pi}\sum_{m=1}^{k-1}\chi (m)D(e^{\frac{2\pi im}{k}}). $$
\end{corollary}
As we mentioned in \Cref{notation df}, for $\chi_{-f}$ the odd quadratic Dirichlet character of conductor $f$,
    Boyd \cite{7} introduced the notation $d_{f}$ for the summation on the right-hand side.
We extend Boyd's notation to all primitive characters by taking  \Cref{warning remarkd-f} into account.
      \begin{notation}\label{dchi}
      Let $\chi$ be an odd Dirichlet character modulo $k$. We define the following notation: $$d_{\chi}\coloneqq\frac{k}{4\pi}\sum_{m=1}^{k-1}\chi (m)D(e^{\frac{2\pi im}{k}}).$$
      \end{notation} 
      From \Cref{coclusion1.2}, if $\chi$ is a primitive odd character, then $d_\chi$ is $L'(\chi,-1)$.
   Furthermore, if $k=f$, where $-f$ is a fundamental discriminant, then $\chi(m)= \chi_{-f}(m)=\big(\frac{-f}{m} \big)$ is the odd quadratic Dirichlet character of conductor $f$. In this case, $d_f$ and $d_{\chi}$ are the same.
      \begin{remark}\label{warning remarkd-f}
          We need to emphasize that the link between $d_{\chi}$ and $L$-functions for primitive characters  $\chi$ mentioned in \Cref{coclusion1.2} is not valid for imprimitive characters. In \Cref{lem 5.2.7}, $d_{\chi}$ for an imprimitive character $\chi$, induced by the primitive character $\chi^*$, is explicitly written in terms of $d_{\chi^*}$, $L'(\chi^*,-1)$, and $L(\bar{\chi},2)$.
      \end{remark}
   Let us first recall a classical fact regarding Dirichlet characters.
 \begin{lemma}[\normalfont \cite{overholt2014course}, Page 91]\label{th 5.4.3}
  Let $\chi$ be a Dirichlet character modulo $k$ with conductor $c$. Then there exists
a unique primitive character $\chi^*$ of conductor $c$ that induces $\chi$. Moreover, $\chi$ is odd if and only if $\chi^*$ is odd.
\end{lemma}

In the following lemma, we announce a linear relation between $d_\chi$ and  $d_{\chi^*}$ with a coefficient in a cyclotomic field depending on the conductor of $\chi^*$. This generalizes \Cref{coclusion1.2} to all characters, providing us with the main tool to prove \Cref{therorm character limit}.
\begin{lemma}[\normalfont \cite{Zagier1991}, Section 9.B]\label{lem 5.2.7}
    Let $\chi$ be a Dirichlet character modulo $k$ and conductor $c$, induced by the primitive Dirichlet character $\chi^*$. Let $q=\frac{k}{c}$ and let $\mu$ be the Möbius function. Define $\gamma$ and  $\beta$ as follows:
\begin{align}\label{gamma}
&\gamma\coloneqq\sum_{d|q}d\,\mu(d) \chi^* (d),
&\beta \coloneqq \frac{ -i \gamma c\tau(\chi^*)}{4\pi\prod_{p \,|\, k}\left(1 - \frac{\overline{\chi^*(p)}}{p^2} \right)}.
\end{align}
 The following equation is valid:
$$d_{\chi}=\gamma d_{\chi^*}=\gamma L'(\chi^*,-1)=\beta L(\overline{\chi},2).$$

\end{lemma}
In the above lemma, the link between $d_\chi$ and $ d_{\chi^*}$ follows from the classical property of the Möbius function, as stated in \cite[Theorem 2.1]{apostol};
    \begin{align*}
\sum_{d|j,d|q}\mu(d)=\sum_{d|\gcd(j,q)}\mu(d)=\begin{cases}
       1 &\quad\text{if} \gcd(j,q)=1\\
       0 &\quad\text{if} \gcd(j,q)>1.
     \end{cases}
  \end{align*} 
Furthermore, if the conductor of $\chi^*$ and the modulo of $\chi$ share the same primes then the constant $\gamma$ in \Cref{lem 5.2.7}
is simply one. This is because we have the equality $\chi(l+mc)=\chi^*(l)\chi_0(l+mc)$, where  $ \chi_0 $ is the principal character modulo $k$ (i.e. $\chi_0(j)=1 \Leftrightarrow \gcd(j,k)=1$). We announce this in the following corollary.
\begin{corollary}\label{specialdchi}
    Let $\chi$ be a character modulo $k$ and conductor $c$,  induced by the primitive Dirichlet character $\chi^*$. Moreover, suppose that $k$ and $c$ have the same prime factors. Then, we have:
$$d_{\chi^*}=d_{\chi}.$$
\end{corollary}

We now proceed with studying $m(P_d)$ and its potential connection to Dirichlet $L$-functions.
In the rest of this section, instead of considering Dirichlet characters modulo $d$ as a function over $\frac{\mathbb{Z}}{d\mathbb{Z}}$, we consider them as functions over the set of $d$-th roots of unity. We explain this process more precisely in \Cref{convendef}, but first, we summarize the important notation and definitions that will be used regularly as follows:
\begin{notation}\label{customary notations}
 Let $d \in \mathbb{Z}_{\geq 1}$ and set the following notation: 
 \begin{itemize}
  \item The set of the $d$-th roots of unity is denoted by $U_d=\{e^{\frac{2\pi i k}{d}}| 0\leq k\leq d-1\}$.
    \item The set of primitive $d$-th roots of unity is denoted by: $$U_d^\times=\{e^{\frac{2\pi i k}{d}}| 1\leq k\leq d-1, \gcd(k,d)=1\}.$$ 
     \item $$
   \frac{\mathbb{Z}}{d\mathbb{Z}}=\{0,1,2,\ldots ,d-1\}.  
$$
\item The group of invertible elements of the ring $ \frac{\mathbb{Z}}{d\mathbb{Z}}$ is denoted by:
\begin{align}\label{U-d*}
  \frac{\mathbb{Z}}{d\mathbb{Z}}^\times=\{k\in \mathbb{Z}| 1\leq k\leq d-1, \gcd(k,d)=1\}. 
 \end{align} 
 \end{itemize}

\end{notation}

It is known \cite[Chapter 7]{stein2011fourier} that the set of Dirichlet characters modulo $d$ form a basis for the $\mathbb{C}$-vector space of functions from $\frac{\mathbb{Z}}{d\mathbb{Z}}^\times$ to $\mathbb{C}$.
Consequently, the set of odd Dirichlet characters modulo $d$ generates $\mathbb{C}$-vector space of functions $f:\frac{\mathbb{Z}}{d\mathbb{Z}}\rightarrow \mathbb{C}$ with the two following properties:
\begin{align*}
&1)f(k)=0 \Leftrightarrow  k\in \frac{\mathbb{Z}}{d\mathbb{Z}}\setminus \frac{\mathbb{Z}}{d\mathbb{Z}}^\times,\\  &2)f(d-x)=-f(x)\ \text{for every}\  x\in  \frac{\mathbb{Z}}{d\mathbb{Z}}.  
\end{align*}
 A Dirichlet character $\chi$ modulo $d$ can be considered as a function $\chi:\frac{\mathbb{Z}}{d\mathbb{Z}} \rightarrow \mathbb{C}$ with some extra properties.
The bijection $\phi_d: \frac{\mathbb{Z}}{d\mathbb{Z}} \rightarrow U_d $, defined by $\phi_d(k)=e^{\frac{2\pi i k}{d}}$ sends $\frac{\mathbb{Z}}{d\mathbb{Z}}^\times$ to the set of primitive $d$-th roots of unity, $U_d^\times$. Instead of considering Dirichlet characters over the set $\frac{\mathbb{Z}}{d\mathbb{Z}}$ we consider them over $\phi_d(\frac{\mathbb{Z}}{d\mathbb{Z}})$ which is $U_d$. We explain our conventional definition more precisely in the following:
\begin{convendef}\label{convendef}
    For $1\leq k\leq d-1$ and a Dirichlet character $\chi$ modulo $d$ we consider the following conventional definition of the Dirichlet character $\chi^\circ\colon U_d \rightarrow\mathbb{C}$, defined by: $$\chi^\circ\left(\left( e^{\frac{2\pi i}{d}}\right)^k\right)\coloneqq \chi\circ\phi_d^{-1} \left(\left( e^{\frac{2\pi i}{d}}\right)^k\right)=\chi(k).$$ 
    \end{convendef}
   Note that the right-hand side of the above equation, $\chi(k)$ is the ordinary definition of Dirichlet characters on $\frac{\mathbb{Z}}{d\mathbb{Z}}$, but the left-hand side $\chi^\circ\left(\left( e^{\frac{2\pi i}{d}}\right)^k\right)$ is the one we consider in the sequel, defined over $U_d$. 
In the following, we state some of the important properties of Dirichlet characters in the new setting which will be needed for the proof of \Cref{therorm character limit}:
\begin{itemize}
    \item Let $\chi^\circ: U_d \rightarrow \mathbb{C},$ we have: 
$$\chi^\circ(z)\neq 0 \Leftrightarrow z\in U_d^\times.$$
\item Let $\chi$ be an odd Dirichlet character modulo $d$, defined over $\frac{\mathbb{Z}}{d\mathbb{Z}}$. For the character $\chi^\circ$ the following properties hold for every $0\leq k\leq d-1$: $$\chi^\circ\left(\overline{e^{\frac{2\pi i k}{d}}}\right)=-\chi^\circ\left(e^{\frac{2\pi i k}{d}}\right).$$
\item For every odd $d$-periodic function $f$, we can also interpret it as a function over $U_d$ via the bijection $\phi_d$ and the process introduced in \Cref{convendef}. For the simplicity, we identify the two functions $f$ over $\frac{\mathbb{Z}}{d\mathbb{Z}}$ and $f^\circ= f \circ \phi_d^{-1}$ over $U_d$  and we use the notation $f$ exclusively, even when considering it over $U_d$.

\item Thanks to the bijection $\phi_d\colon \frac{\mathbb{Z}}{d\mathbb{Z}} \rightarrow U_d$ we conclude that the set of  
  Dirichlet characters modulo $d$ over $U_d$ forms a basis for the $\mathbb{C}$-vector space $\mathbb{C}^{ U_d^\times}$. 
 Furthermore, the set of odd Dirichlet characters modulo $d$ over $U_d$ generates the set of functions $f:U_d\rightarrow \mathbb{C}$ with the two following properties:
\begin{align}\label{genrarting-characters}
&1)f(z)=0 \Leftrightarrow  z\in U_d\setminus U_d^\times,\\  &2)f(\bar{z})=-f(z)\ \text{for every}\  z\in  U_d.\nonumber  \end{align}
\end{itemize}

\begin{definition}\label{oddmoddcomplex}
  A function $f:U_d\rightarrow \mathbb{C}$ is called \enquote{odd} if it satisfies 
  $f(\bar{z})=-f(z)$ for every $z$ in $U_d$. In our context, a function is called \enquote{even} if it is invariant under complex conjugation.
\end{definition}
Using the above terminology, the set of odd Dirichlet characters modulo $d$ over $U_d$ generates the set of odd functions whose support is $U_d^\times$. Let $d \in \mathbb{N}$ and $k$ be a divisor of $d$. In the following lemma, we construct an odd function on $U_k$ with support $U_k^\times$ from an odd function defined on $U_d$.

\begin{proposition}\label{lem 5.4.1}
Let $d \in \mathbb{N}$ and $k$ be a divisor of $d$. Let $f_d:U_d \rightarrow \mathbb{Z}$ be an odd function and let  $1_{U_k^\times}$ be the characteristic function of $V_k$. Let $\hat{f}_{k,d}:U_d \rightarrow \mathbb{Z}$ be defined by $\hat{f}_{k,d}\coloneqq1_{U_k^\times}f_d$. Let  $f_{k,d}$ be the restriction of $\hat{f}_{k,d}$ to $U_k$. Then $f_{k,d}$ is an odd function and is written uniquely in terms of odd Dirichlet characters modulo $k$ over $U_k$.
\end{proposition}
\begin{Proof}
The characteristic function $1_{U_k^\times}$ is even and $f_d$ is odd, therefore the product is odd. Thus it is written uniquely in terms of odd Dirichlet characters modulo $k$ over $U_k$.
\end{Proof}

We now come back to our problem of writing $m(P_d)$ in terms of $L$-functions. In \Cref{eq 4.1.11} the Mahler measure of $P_d$ is expressed in terms of $S_d$.  Let us prove the following proposition concerning $S_d$: 

\begin{proposition}\label{th 5.0.2}
Let $d\in \mathbb{Z}_{\geq 1}$, and $S_d=3\sum_{1 \leq k \leq d-1}(d-2k)D\left(e^{\frac{2\pi i k}{d}}\right)$. Then, for every primitive odd Dirichlet character $\chi$ of conductor $k$, such that $k|d$, there exists a coefficient $C_{k,\chi}^d$ in $\mathbb{Q}\left(e^{\frac{2\pi i}{\phi(d)}}\right)$ such that:
$$S_d=\sum_{k|d}\sum_{\chi\,\text{primitive odd mod $k$} } C_{k,\chi}^d L'(\chi,-1).$$
\end{proposition}
\begin{Proof}
 We consider $S_d$ as an inner product of two vectors in $\mathbb{C}^{d-1} $ as follows:
$$S_d=3\left<\left[d-2, \ldots, 2-d\right], \left[D\left(e^{\frac{2\pi i}{d}}\right), \ldots , D\left(e^{\frac{2\pi i(d-1)}{d}}\right)\right] \right>.$$
In \Cref{coclusion1.2} and \Cref{lem 5.2.7} we have seen the link between $d_\chi$ and $L$-functions for primitive odd and imprimitive odd characters, respectively. Thus, writing $S_d$ in terms of a linear combination of $d_\chi$, for some odd $\chi$ leads to a representation of $S_d$ in terms of $L$-functions. We write $d_\chi$ as follows:
$$d_{\chi}=\frac{d}{4\pi}\left<\left[\chi(1), \ldots, \chi(d-1)\right], \left[D\left(e^{\frac{2\pi i}{d}}\right), \ldots , D\left(e^{\frac{2\pi i(d-1)}{d}}\right)\right] \right>.$$
Thanks to the bilinearity of the inner product, it suffices to write the coefficient vector $$C_d \coloneqq \left[d-2, \ldots, 2-d\right].$$ as a linear combination of character vectors, which are vectors of the form  

 $$\boldsymbol{\chi_d}
\coloneqq \left[\chi(1), \ldots, \chi(d-1)\right],$$  for some odd Dirichlet character $\chi$. 
By using the bijection $\phi_d\colon \frac{\mathbb{Z}}{d\mathbb{Z}} \rightarrow U_d$, we write 
 character vectors in the following form:  

 $$\boldsymbol{\chi^\circ_d}
\coloneqq \left[\chi^\circ\left(e^{\frac{2\pi i}{d}}\right), \ldots, \chi^\circ\left(e^{\frac{2\pi i(d-1)}{d}}\right)\right],$$  for the odd $\chi^\circ$, constructed from $\chi$ using \Cref{convendef}. To use linear algebra techniques we define the following function:
 \begin{alignat*}{3}
 f_d\colon &U_d &&\longrightarrow \quad&&\mathbb{Z}\\
&e^{\frac{2\pi i k}{d}}&&\longmapsto &&d-2k,
    \end{alignat*}
where $1\leq k\leq d-1$. Using the above function, the question of writing $C_d$ as a linear combination of vectors $\boldsymbol{\chi^\circ_d}$ is equivalent to writing the function $f_d$ as 
a linear combination of odd Dirichlet characters over $U_d$.
We note that $f_d$ is an odd function, because for any $z=e^{\frac{2k\pi i}{d}}\in U_d$ the following holds:
$$f_d(\bar{z})=f_d\left(e^{\frac{2\pi i(d-k) }{d}}\right)=d-2(d-k)=2k-d=-f_d\left(e^{\frac{2\pi i k}{d}}\right)=-f_d(z).$$
However, $f_d$ does not satisfy the first condition of \Cref{genrarting-characters}. In other words we have $f_d(z) \neq 0$ even for $z\in U_d\setminus U_d^\times$. Thus we can not write $f_d$ only as a linear combination of odd Dirichlet characters modulo $d$.
To remove this difficulty, we write $f_d=\sum_{k|d}\hat{f}_{d,k}$, with $\hat{f}_{d,k}=1_{V_k}f_d$, as before. We notice that $U_d=\bigcup_{k|d}V_k$, where the union is disjoint. One can restrict each $\hat{f}_{d,k}$ to $U_k$ by considering $f_{d,k}:U_k \rightarrow \mathbb{C}$ in the same way as introduced in \Cref{lem 5.4.1}.
Therefore, thanks to \Cref{lem 5.4.1} for each $k|d$, $f_{d,k}$ can be represented uniquely as a linear combination of odd Dirichlet characters $\chi^\circ$ modulo $k$ over $U_k$. Finally, using the bijection $\phi_k^{-1}:U_k  \rightarrow \frac{\mathbb{Z}}{k\mathbb{Z}}$, we can write $S_d=\sum_{k|d}\sum_{\chi \text{odd of mod $k$}}C_{k,\chi}^dd_\chi$. To explain that $C_{k,\chi}^d\in \mathbb{Q}\left(e^{\frac{2\pi i}{\phi(d)}}\right)$,
we recall that the vector $[d-2, \dots, 2-d]\in \mathbb{Z}^d$  is written uniquely as a linear combination of Dirichlet characters modulo $k$, with $k|d$. On the other hand, the values of Dirichlet characters modulo $k$ belong to the cyclotomic field $\mathbb{Q}\left(e^{\frac{2\pi i}{\phi(k)}}\right)$. Thanks to the multiplicativity of Euler's function (see \cite[Page 13]{deza2021mersenne}), for any divisor $k$ of $d$ we have $\phi(k)|\phi(d)$, so $\mathbb{Q}\left(e^{\frac{2\pi i}{\phi(k)}}\right) \subseteq \mathbb{Q}\left(e^{\frac{2\pi i}{\phi(d)}}\right)$. Using the fact that Dirichlet characters modulo $d$ construct a basis for $\mathbb{C}^{U_d^\times}$ we conclude that all $C_{k,\chi}^d$ belonging to $\mathbb{Q}\left(e^{\frac{2\pi i}{\phi(d)}}\right)$.
Using \Cref{lem 5.2.7}, for a character $\chi$ modulo $k$ and conductor $c$ induced by $\chi^*$ we can write $d_\chi=\gamma_\chi L'(\chi^*,-1)$, where $\gamma_\chi$ is constructed from $\chi$ according to \Cref{gamma} and we have: $$\gamma_\chi \in \mathbb{Q}\left(e^{\frac{2\pi i}{\phi(c)}}\right)\subseteq \mathbb{Q}\left(e^{\frac{2\pi i}{\phi(k)}}\right) \subseteq \mathbb{Q}\left(e^{\frac{2\pi i}{\phi(d)}}\right).$$
This implies that $S_d$ can be written as a linear combination of $L$-functions with coefficients in the cyclotomic field $\mathbb{Q}\left(e^{\frac{2\pi i}{\phi(d)}}\right)$.
\end{Proof}
The above proposition is the key to proving \Cref{therorm character limit}.
\begin{Proof}{Proof of \Cref{therorm character limit}:}
From \Cref{eq 4.1.11} and 
 \Cref{th 5.0.2} we can write $S_{d+1}$ and $S_{d+2}$ in terms of $L$-functions associated with primitive odd characters of conductor $k$ where $k|(d+1)$ or $k|(d+2)$, respectively. This completes the proof.
\end{Proof}
 
 \begin{remark}\label{nonunicrep}
     We notice that the representation of $m(P_d)$ in terms of $L$-functions whose existence is guaranteed by \Cref{therorm character limit} is not unique. One can uniquely express the coefficient vector $C_d= [d-2, \dots,  2-d]$ as a linear combination of primitive odd Dirichlet characters. However, in the formula of $S_d$ we compute the inner product of $C_d$ with the dilogarithm vector $$D_d\coloneqq\left[D\left(e^{\frac{2\pi i}{d}}\right), \dots , D\left(e^{\frac{2\pi i(d-1)}{d}}\right)\right],$$ which changes the situation. More precisely, for any vector $V \in D_d^\bot$, perpendicular to the vector $D_d$, we have $S_d= 3 \langle C_d+V, D_d\rangle=3 \langle C_d,D_d\rangle$.
  \end{remark}
\section{Solutions to Chinburg's conjectures using the \texorpdfstring{$P_d$}{} family}\label{Solutions to Chinburg's conjecture-sec}
As we have seen in the previous section, \Cref{therorm character limit} revealed the links between  $m(P_d)$ and $L$-functions of primitive odd Dirichlet characters. \Cref{6.3.6} and \Cref{genchinsolution} shows the solutions that we constructed using $P_d$ family. In this section, we explain the method to find these solutions. \\

The Mahler measure of $P_d$ thanks to \Cref{eq 4.1.11} is written in terms of $S_d$. According to \Cref{th 5.0.2}, for a fixed $t\in \mathbb{N}$ we can express $S_t$ in terms of $L$-functions of primitive characters of conductor $f$, where $f|t$. Furthermore, an inductive process can express $S_t$ as a linear combination of the Mahler measure of $P_d$ polynomials, with $d\leq t$.
 Suppose there exists only one primitive odd Dirichlet character of conductor $f$, namely $\chi_{-f}$, then by considering a suitable linear combination of $S_d$, where $d\leq t$, one can express $L'(\chi_{-f},-1)$ in terms of an integral linear combination of the Mahler measures of $P_d$ polynomials with $d\leq t$. Thanks to the property of the Mahler measure that $m(PQ)=m(P)+m(Q)$ and $m(P/Q)=m(P)-m(Q)$, for $Q\neq 0$, an integral linear combination of Mahler measures of members of $P_d$ family is itself the Mahler measure of a rational function involving product and division of $P_d$ polynomials, with $d\leq t$. 
 However, if there is more than one primitive Dirichlet character of conductor $f$ then we may not be able to propose any solution to \Cref{chinberg weak} using $P_d$ polynomials. It might still be possible to propose a solution to \Cref{genchin}. More precisely, using a similar argument if there exists only one pair of complex conjugate primitive odd Dirichlet character of conductor $f$, namely $\chi$ and $\bar{\chi}$, $P_d$ can propose a solution to \Cref{genchin} for $\chi$.\\

Our computation of $S_d$ in terms of $L$-functions involves applying \Cref{coclusion1.2}, \Cref{lem 5.2.7}, and linear algebra techniques. We detail the computation for $S_{20}$
  as an example, as the others follow a similar process. For the reader's convenience, a table presenting the 
$S_d$ in terms of $L$-functions for the values of 
$d$ needed for our computation, is provided in \Cref{appendix}.

\begin{proposition}\label{S20-S24prop}
   The formula presented in  \Cref{eq35} in \Cref{table-of-sd} holds.
\end{proposition}
\begin{Proof}

\begin{align*}
    \frac{S_{20}}{2\pi}=&\frac{3}{\pi}\left(18D\left(e^{\frac{2\pi i}{20}}\right)+16D\left(e^{\frac{4\pi i}{20}}\right)+14D\left(e^{\frac{6\pi i}{20}}\right)+12D\left(e^{\frac{8\pi i}{20}}\right)\right)+\\
    &+\frac{3}{\pi}\left(10D\left(e^{\frac{10\pi i}{20}}\right)+8D\left(e^{\frac{12\pi i}{20}}\right)+6D\left(e^{\frac{14\pi i}{20}}\right)+4D\left(e^{\frac{16\pi i}{20}}\right)+2D\left(e^{\frac{18\pi i}{20}}\right)\right).
\end{align*}
There is only one primitive odd character mod $20$, the quadratic one. We write $S_{20}$ in terms of $d_4$, $d_{\chi_5(2,.)}$, $d_{\chi_5(3,.)}$ and $d_{20}$. By using \Cref{s5} and \Cref{s4} we have:
\begin{align*}
    \frac{S_{20}}{2\pi}=&\re\left(\frac{72-24i}{5}L'({\chi_5(2,.)},-1)\right)+15L'(\chi_{-4},-1)\\
    &+\frac{3}{\pi}\left(16D\left(e^{\frac{2\pi i}{10}}\right)+8D\left(e^{\frac{6\pi i}{10}}\right)\right)\\
    &+\frac{3}{\pi}\left(18D\left(e^{\frac{2\pi i}{20}}\right)+14D\left(e^{\frac{6\pi i}{20}}\right)+6D\left(e^{\frac{14\pi i}{20}}\right)+2D\left(e^{\frac{18\pi i}{20}}\right)\right).
\end{align*}
To write $\frac{3}{\pi}\left(16D\left(e^{\frac{2\pi i}{10}}\right)+8D\left(e^{\frac{6\pi i}{10}}\right)\right)$ in terms of $d_{\chi}$, we need to use all the odd characters modulo $10$ which are the imprimitive characters $\chi_{10}(3,.)$ and $\chi_{10}(7,.)$. $$d_{\chi_{10}(7,.)}=\overline{d_{\chi_{10}(3,.)}}=\frac{5}{\pi}\left(D\left(e^{\frac{2\pi i}{10}}\right)-iD\left(e^{\frac{6\pi i}{10}}\right)\right).$$
The characters $\chi_{10}(3,.)$ and $\chi_{10}(7,.)$  are respectively induced by $\chi_{5}(3,.)$ and $\chi_{5}(2,.)$. An application of \Cref{lem 5.2.7} implies that $d_{\chi_{10}(3,.)}=(1+2i)d_{\chi_{5}(3,.)},$ and  $d_{\chi_{10}(7,.)}=(1-2i)d_{\chi_{5}(2,.)}$. Thus, we have the following identity: 

\begin{align*}
 \frac{3}{\pi}\left(16D\left(e^{\frac{2\pi i}{10}}\right)+8D\left(e^{\frac{6\pi i}{10}}\right)\right)&=\frac{48-24i}{10}d_{\chi_{10}(3,.)}+\frac{48+24i}{20}d_{\chi_{10}(7,.)}\\
 &=\frac{48+36i}{5}d_{\chi_{5}(3,.)}+\frac{48-36i}{5}d_{\chi_{5}(2,.)}
  \\
 &=\re\left(\frac{96-72i}{5}L'({\chi_{5}(2,.),-1)}\right).
\end{align*}
We write $\frac{3}{\pi}\left(18D\left(e^{\frac{2\pi i}{20}}\right)+14D\left(e^{\frac{6\pi i}{20}}\right)+6D\left(e^{\frac{14\pi i}{20}}\right)+2D\left(e^{\frac{18\pi i}{20}}\right)\right)$ in terms of $d_{\chi}$ associated with all the odd characters modulo $20$. The odd characters modulo $20$ are $\chi_{-20}$ which is primitive and $\chi_{20}(11,.)$, 
$\chi_{20}(13,.)$, and $\chi_{20}(17,.)$ which are induced by $\chi_{-4}$, $\chi_{5}(3,.)$ and $\chi_{5}(2,.)$, respectively:
\begin{align*}
&d_{20}=\frac{10}{\pi}\left(D\left(e^{\frac{2\pi i}{20}}\right)+D\left(e^{\frac{6\pi i}{20}}\right)+D\left(e^{\frac{14\pi i}{20}}\right)+D\left(e^{\frac{18\pi i}{20}}\right)\right),\\
 &d_{\chi_{20}(13,.)}=\overline{d_{\chi_{20}(17,.)}}=\frac{10}{\pi}\left(D\left(e^{\frac{2\pi i}{20}}\right)+iD\left(e^{\frac{6\pi i}{20}}\right)-iD\left(e^{\frac{14\pi i}{20}}\right)-D\left(e^{\frac{18\pi i}{20}}\right)\right),\\
 &d_{\chi_{20}(11,.)}=\frac{10}{\pi}\left(D\left(e^{\frac{2\pi i}{20}}\right)-D\left(e^{\frac{6\pi i}{20}}\right)-D\left(e^{\frac{14\pi i}{20}}\right)+D\left(e^{\frac{18\pi i}{20}}\right)\right).
  \end{align*}
  By solving a system of linear equations we have:
\begin{align*}
 &\frac{3}{\pi}\left(18D\left(e^{\frac{2\pi i}{20}}\right)+14D\left(e^{\frac{6\pi i}{20}}\right)+6D\left(e^{\frac{14\pi i}{20}}\right)+2D\left(e^{\frac{18\pi i}{20}}\right)\right)=\\
 &3d_{20}+\frac{6-3i}{5}d_{\chi_{20}(13,.)}+\frac{6+3i}{5}d_{\chi_{20}(17,.)}.   
\end{align*}
An application of \Cref{lem 5.2.7} gives
$ d_{\chi_{20}(13,.)}=(1+2i)d_{\chi_{5}(3,.)},$ and $d_{\chi_{20}(17,.)}=(1-2i)d_{\chi_{5}(2,.)}.$ 
Finally, after simplification, we obtain \Cref{eq35}.
\end{Proof}

\section{Perspective}\label{pers}

\Cref{6.3.6} convinces us that $P_d$ has the potential to produce further solutions to various versions of Chinburg's conjectures. 
The direct computation of $m(P_d)$ has become harder when $d$ increases. Moreover, the representation is not unique as we explained in \Cref{nonunicrep}.  
 Hence, one important future project is implementing an algorithm that generates these solutions for different versions of Chinburg's conjectures for the case $n=1$.\\

 The key property of $P_d$ that makes the bridge between $m(P_d)$ and Dirichlet $L$-functions is that $m(P_d)$ is written in terms of the dilogarithm function evaluated at roots of unity which allows us to apply \Cref{df} or \Cref{coclusion1.2}. This is thanks to the exactness of $P_d$ and that its toric points (roots of the polynomial located on the complex 2-dimensional Torus) are roots of unity (see \cite{mehrabdollahei2021mahler} for more details).  Hence, one can hope to construct other families of exact polynomials with this property and maybe combine these polynomials to produce even more solutions to Chinburg's conjectures. Other such families already have been constructed by Ray\cite{Ray1987RelationsBM}, Boyd, and Rodriguez-Villegas\cite{boyd_rodriguez-villegas_2002} and have been indeed applied fruitfully to Chinburg's conjecture as mentioned in the introduction. A work in progress of David Hokken, the second author of the present paper, and Berend Ringeling provides more insight into Chinburg's conjectures by considering the mentioned perspectives.\\
 Furthermore, as mentioned in \cref{mostgeneralchin}, one can investigate the most general case of Chinburg's conjecture for all characters and construct examples using the $P_d$
  family, known polynomial families from the literature, and Ray's work \cite[Proposition 22]{Ray1987RelationsBM}.

\section{Appendix}\label{appendix}
\subsection{Conrey label of Dirichlet characters used in LMFDB}
We use the notation $\chi_{q}(n,\cdot)$ to identify Dirichlet characters $\mathbb{Z}\to \mathbb{C}$, where $q$ is the 
modulus, and $n$ is the index, a positive integer coprime to $q$ that identifies a Dirichlet character modulo $q$ as described below.  The LMFDB label $\texttt{q.n}$, with $1\le n < \max(q,2)$ uniquely identifies $\chi_{q}(n,\cdot)$.
Introduced by Brian Conrey, this labeling system is based on an explicit isomorphism between the multiplicative group $(\mathbb{Z}/q\mathbb{Z})^\times$ and the group of Dirichlet characters modulo $q$ that makes it easy to recover the
order, the
conductor, and the
parity
of a Dirichlet character from its label, or to induce characters.
As an example,
$\chi_q(1, \cdot)$ is always trivial, $\chi_q(m,\cdot)$ is real if $m^2=1\bmod q$, and for all $m,n$ coprime to $q$ we have $\chi_q(m,n)=\chi_q(n,m)$.
\begin{description}
    \item[Modulo Prime Powers:] 
For prime powers $q=p^e$ we define $\chi_q(n,\cdot)$ as follows.

\begin{itemize}

\item For each odd prime $p$ we choose the least positive integer $g_p$ which
 is a
  primitive root for all $p^e$, and then for $n \equiv g_p^a $ mod $p^{e}$ and $m
  \equiv g_p^{b} $ mod $p^{e}$ coprime to $p$ we define
  $$
    \chi_{p^e}(n, m) = \exp\left(2\pi i \frac{a b}{\phi(p^{e})} \right).
  $$

\item We note that $\chi_2(1, \cdot)$ is the trivial character modulo $2$, $\chi_4(3, \cdot)$ is the
  unique nontrivial character modulo $4$, and for larger powers of $2$ we choose
  $-1$ and $5$ as generators of the multiplicative group. For $e > 2$, if
  $$
    n \equiv \epsilon_a 5^a \pmod{2^e}
  $$
  and
  $$
    m \equiv \epsilon_b 5^b \pmod{2^e}
  $$
  with $\epsilon_a, \epsilon_b \in \{\pm 1\}$, then
  $$
    \chi_{2^e}(n, m) = \exp\left(2 \pi i \left(\frac{(1 - \epsilon_a)(1 - \epsilon_b)}{8}
        + \frac{ab}{2^{e-2}}\right)\right).
  $$

\end{itemize}
 \item[Modulo General Positive   Integers:] For general $q$, the function $\chi_q(n, m)$ is defined multiplicatively: $$\chi_{q_1 q_2}(n, m) := \chi_{q_1}(n, m)\chi_{q_2}(n, m)$$  for all coprime positive integers $q_1$ and $q_2$.  The Chinese remainder theorem implies that this definition is well-founded and that every Dirichlet character can be defined in this way.  In particular, every Dirichlet character $\chi$ modulo $q$ can be written uniquely as a product of Dirichlet characters of prime power modulus.
\end{description}

\subsection{Tables}

The constants $C_{k,n}^d$, associated with the coefficient of $L'(\chi_{k}(n,.),-1)$, appearing in \Cref{tablemp2}  in the introduction as well as the following table are introduced in \Cref{table-of-constant}.\vspace{.3cm}

\small{
\begin{center}
    \parbox{1.019\textwidth}{%
\fbox{%
    \parbox{1.2\textwidth}{%
    {\small
    \begin{align}
  &\frac{S_3}{2\pi}= 2L'(\chi_{\chi_{-3}},-1) \label{s3}\\
&\frac{S_4}{2\pi}= 3L'(\chi_{\chi_{-4}},-1)\label{s4}\\
& \frac{S_5}{2\pi} =\re\left(\frac{18-6i}{5}L'(\chi_5(2,.),-1)\right)\label{s5}\\
 &\frac{S_6}{2\pi}=16 L'(\chi_{-3},-1)\label{s6}\\
   & \frac{S_7}{2\pi}=2L'(\chi_{-7},-1) + \re\left(\frac{16-8\sqrt{3}i}{7}L'(\chi_{7}(3,.),-1)\right) \label{s7}\\
     &\frac{S_8}{2\pi}=\frac{15}{2}L'(\chi_{-4},-1)+3L'(\chi_{-8},-1)\label{s8}\\
    &\frac{S_9}{2\pi}=\re\left(\frac{12-4\sqrt{3}i}{3}L'(\chi_9(2,.),-1)\right)+\frac{20}{3}L'(\chi_{-3},-1)\label{s9}\\
     &\frac{S_{10}}{2 \pi}=\re\left(\frac{84-48i}{5}L'({\chi_5(2,.)},-1)\right)\label{s10}\\
       &\frac{S_{11}}{2\pi}=2\re\left(C_{11,2}^9L'({\chi_{11}(2,.)},-1)+C_{11,7}^9L'({\chi_{11}(7,.)},-1)\right)+\frac{6}{5}L'({\chi_{-11}},-1)\label{s11}\\
       & \frac{S_{12}}{2\pi}=21L'(\chi_{-4},-1)+38L'(\chi_{-3},-1)\label{s12}\\
        &\frac{S_{13}}{2\pi}=2\re\left( C_{13,5}^{11}L'({\chi_{13}(5,.)},-1)+C_{13,2}^{11}L'(\chi_{13}(2,.),-1)+C_{13,
        6}^{11}L'(\chi_{13}(6,.),-1)\right)\label{S13}\\ 
        & \frac{S_{14}}{2\pi}=\re\left(\frac{92-60 \sqrt{3}i}{7} L'({\chi_{7}(3,.)},-1)\right)+ 4L'(\chi_{-7},-1)\label{S14}\\
        &\frac{S_{15}}{2\pi}=\re\left(\frac{96-12i}{5}L'(\chi_5(2,.),-1)\right)+3L'(\chi_{-15},-1)+16L'(\chi_{-3},-1)\label{s15}\\
&\frac{S_{16}}{2\pi}=\re\left((3+3i)L'(\chi_{16}(3,.),-1)\right)+\frac{15}{2}L'({\chi_{-8}},-1)+\frac{63}{4}L'({\chi_{-4}},-1)\label{S16}\\
& \frac{S_{17}}{2\pi}=2\re\left(C_{17,3}^{15}L'({\chi_{17}(3,.)},-1)+C_{17,5}^{15}L'({\chi_{17}(5,.)},-1)+C_{17,10}^{15}L'({\chi_{17}(10,.)},-1)+C_{17,12}^{15}L'({\chi_{17}(12,.)},-1)\right)\label{s17}\\
    &\frac{S_{18}}{2\pi}=\re\left(\frac{36-20\sqrt{3}i}{3}L'(\chi_9(2,.),-1)\right)+\frac{160}{3}L'(\chi_{-3},-1)\label{s18}\\
    &\frac{S_{20}}{2\pi}=3L'(\chi_{20},-1)+ \re\left(\frac{192-114i}{5}L'({\chi_{5}(2,.)},-1)\right)+15L'(\chi_{-4},-1),\label{eq35}\\
        &\frac{S_{24}}{2\pi}=3L'(\chi_{-24},-1)+9L'(\chi_{-8},-1)+ \frac{105}{2}L'(\chi_{-4},-1) +79L'(\chi_{-3},-1).\label{eq36}
        \end{align}} }}
}\captionsetup{type=figure}
\captionof{table}{Table of  $S_d$ in terms of $L$-functions.}
\label{table-of-sd}
\end{center}}

\begin{center}
    \parbox{1.019\textwidth}{%
\fbox{%
    \parbox{\textwidth}{%
    {\small
    \begin{align*}
  &C_{11,2}^9 = C_{11,2}^{10}=\frac{51+21\sqrt{5}}{55} + i\frac{(3\sqrt{5}-21)\sqrt{10-2 \sqrt{5}}}{110} \nonumber\\
&\\
&C_{11,7}^9=C_{11,7}^{10}= \frac{51-21\sqrt{5}}{55} + 
i\frac{(6\sqrt{5}+9)\sqrt{10-2 \sqrt{5}}}{55} \\
&\\
        &C_{13,5}^{11}=C_{13,5}^{12}=1+i\nonumber\\
    &\\
    &C_{13,2}^{11}=C_{13,2}^{12}=\frac{10+4\sqrt{3}}{13}- i\frac{6\sqrt{3}+2}{13} \nonumber\\
&\\
        &C_{13,6}^{11}=C_{13,6}^{12}=\frac{10-4\sqrt{3}}{13}+i\frac{6\sqrt{3}-2}{13}\\
        & \\
 &C_{17,3}^{15}=C_{17,3}^{16}=\frac{45+18\sqrt{2}-(15\sqrt{2}-39)\sqrt{\sqrt{2} + 2}}{68}+i\frac{27+21\sqrt{2}-(6\sqrt{2}+15)\sqrt{2-\sqrt{2}}}{68}\\
 &\\
&C_{17,5}^{15}=C_{17,5}^{16}=\frac{45-18\sqrt{2}-(24\sqrt{2}-9)\sqrt{\sqrt{2} + 2}}{68}+i\frac{27-21\sqrt{2}-(9\sqrt{2}+3)\sqrt{2-\sqrt{2}}}{68}\\
&\\
&C_{17,10}^{15}=C_{17,10}^{16}=\frac{45-18\sqrt{2}+(24\sqrt{2}-9)\sqrt{\sqrt{2} + 2}}{68}+i\frac{-27+21\sqrt{2}-(9\sqrt{2}+3)\sqrt{2-\sqrt{2}}}{68}\\
&\\
&C_{17,12}^{15}=C_{17,12}^{16}=\frac{45+18\sqrt{2}+(15\sqrt{2}-39)\sqrt{\sqrt{2} + 2}}{68}+i\frac{-27-21\sqrt{2}-(6\sqrt{2}+15)\sqrt{2-\sqrt{2}}}{68}
        \end{align*}} }}
}\captionsetup{type=figure}
\captionof{table}{Table of the constants appearing in \Cref{tablemp2}.}
\label{table-of-constant}
\end{center}
\vspace{.7cm}
The following two tables present partial solutions to Chinburg's conjecture for the known conductors.

\begin{center}
    \parbox{1.019\textwidth}{%
\fbox{%
    \parbox{\textwidth}{%
    {\small
    
\begin{align*}
&d_3=m(x+y+1)
 \\
 &\\
&d_4=\frac{1}{2}m\left((x+1)^2y^2+(x-1)^2\right)
  \\
   &\\
&d_7=\frac{7}{8}m\left(\frac{(x^7-1)}{x-1}(y-1)^2+7x^2(x+1)^2y\right)
 &\\
  \\
&d_8=m\left((x^4+1)(y-1)^2+8x^2y\right)\\ 
 &\\
&d_{20}=\frac{5}{2}m\left((x^8-x^6+x^4-x^2+1)(y-1)^2+20x^2(x^2-1)^2y\right)\\
 &\\
&d_{24}=3m\left((x^8-x^4+1)(y-1)^2+24x^2(x^2-1)^2y\right)
 \end{align*}} }}
}
\captionsetup{type=figure}
\captionof{table}{Table of the solutions to Chinburg's conjecture computed by Ray in \cite{Ray1987RelationsBM}.}\label{Table of Rayexamples}
\end{center}

\begin{center}

    \parbox{1.019\textwidth}{%
\fbox{%
    \parbox{\textwidth}{%
    {\small
    \begin{align*}
&d_{11}\sim \frac{3}{2}m((x+1)^2(x^2+x+1)y-(x^2-x+1)^2)
 \\
   &\\
&d_{15}\sim 6m((x+1)^2y-(x^2-x+1))
  \\
    &\\
&d_{19}\sim \frac{5}{2}m((x^5+x^4+x^3+x^2+x)y^2-(x^6+6x^5+2x^4-8x^3+2x^2+6x+1)y+ (x^5+x^4+x^3+x^2+x)) 
 \\
   &\\
 &d_{23}\sim 6 m((x^4-x^2+1)y^2+(x^6-6x^5+12x^3-6x+1)y+(x^4-x^2+1)x^2) 
 \\
   &\\
&d_{35}\sim 12m((x^2+x+1)^3y-(x^2+1)(x^2-x+1))
 \\
   &\\
&d_{39}\sim 18m((x+1)^4y-(x^2+x+1)(x^2-x+1))
 \\
   &\\
&d_{40}\sim 6m((x^2-x+1)(x^2+1)y^2+x(14x^2-32x+14)y+(x^2-x+1)(x^2+1))
  \\
    &\\
 &d_{55}\sim 30m((x+1)^2(x^2+x+1)y-(x^4-x^3+x^2-x+1))
  \\
    &\\
 &d_{84}\sim 36m((x+1)^4y-(x^2+1)(x^2-x+1)) 
  \\
    &\\
 &d_{120}\sim 36m((x^2+1)(x^2+x+1)y^2+2(x^2-3x+1)(x^2+4x+1)y+(x^3+1)(x^2+x+1))
  \\
    &\\
 &d_{303}\sim 132 m((x^8+x^7+x^6+x^2+x+1)(y^2+1)+(2x^8+2x^7-98x^4+2x^3-49x^2+2x+2)y)
  \\
  &\\
 &d_{755}\sim 410m((x^8+x^6+x^4+x^2+1)(y^2+1)+(2x^8-37x^6+5x^5+70x^4+5x^3-37x^2+2)y)
 \end{align*}} }}
}

\captionsetup{type=figure}
\captionof{table}{Table of solutions to Chinburg's conjecture with approximated constants \cite{boyd_rodriguez-villegas_2002,Mahler2,Luo_2019}. }
\label{Numericalyverifiedtable}

\end{center}

\bibliographystyle{acm}
\bibliography{thebib}
\section*{}
Marie José Bertin, Department of Mathematics at IMJ-PRG (Institut de Mathématiques de Jussieu-Paris Rive Gauche), Sorbonne University, 4 Place Jussieu, 75252 Paris, CEDEX 5, France.\\
\textit{Email address:}
{\color{blue}marie-jose.bertin@imj-prg.fr}\\
\textit{ORCID Number:https://orcid.org/0000-0002-9331-2737}\\

\noindent Mahya Mehrabdollahei, Department of Mathematics,
Georg-August-Universität Göttingen,
Bunsenstraße 3-5,
D-37073 Göttingen,
Germany.\\
\textit{Email address:}
{\color{blue}mahya.mehrabdollahei@mathematik.uni-goettingen.de}\\
\textit{ORCID Number:}https://orcid.org/0000-0002-0323-0405

\end{document}